\documentclass[12pt,reqno]{article}
\usepackage{dcolumn}
\usepackage{amsthm,epsfig}
\usepackage{array,dcolumn,graphicx,amsmath}
\usepackage{booktabs,lscape,multirow,psfrag,rotating}
\usepackage{amssymb,amsmath,amsthm}
\usepackage{longtable}
\usepackage{booktabs,lscape,multirow,rotating}
\usepackage[normal]{caption2}
\usepackage[mathscr]{eucal}
\usepackage{eucal}
\usepackage{amssymb,pstricks}
\newcommand{\mXbar}{\overline{\bf{X}}}

\newcommand{\bSigma}{\boldmath$\Sigma$}
\newcommand{\mbSigma}{\mbox{\bSigma}}
\newcommand{\bTheta}{\boldmath$\Theta$}
\newcommand{\mbTheta}{\mbox{\bTheta}}

\newcommand{\bGamma}{\boldmath$\Gamma$}
\newcommand{\mbGamma}{\mbox{\bGamma}}
\newcommand{\bdelta}{\boldmath$\delta$}
\newcommand{\mbdelta}{\mbox{\bdelta}}

\newcommand{\bmu}{\boldmath$\mu$}
\newcommand{\mbmu}{\mbox{\bmu}}

\newcommand{\bpi}{\boldmath$\pi$}
\newcommand{\mbpi}{\mbox{\bpi}}

\newcommand{\bgamma}{\boldmath$\gamma$}
\newcommand{\mbgamma}{\mbox{\bgamma}}

\newcommand{\bPsi}{\boldmath$\Psi$}
\newcommand{\mbPsi}{\mbox{\bPsi}}

\newcolumntype{d}[1]{D{.}{.}{#1}}
\begin{document}
\def\s#1{\oalign{$#1$\crcr\hidewidth \normal$\sim$ \hidewidth}}
\thispagestyle{empty}
\newtheorem{theorem}{\indent THEOREM}
\newtheorem{prop}{\indent Proposition}
\newtheorem{lemma}{\indent Lemma}
\renewcommand{\proofname}{\hspace*{\parindent}{Proof.}}

\begin{center}
{\bf  On the consistent estimators of the population covariance matrix and its reparameterizations}
\end{center}

\begin{center}
{Ming-Tien Tsai and Chia-Hsuan Tsai}

\vspace{0.2cm}
{\it Institute of Statistical Science, Academia Sinica, Taipei, Taiwan 11529, R.O.C.} \\

\end{center}

\vspace{0.3cm}
\begin{center}
\parbox{14cm}{\sloppy \, \,
{\bf Abstract}. For the high-dimensional covariance estimation problem, when $\lim_{n\to \infty}p/n=c \in (0,1)$ the
orthogonally equivariant estimator of the population covariance matrix proposed by Tsai and Tsai (2024b) enjoys some
optimal properties. Under some regularity conditions, they showed that their novel estimators of eigenvalues are consistent
with the eigenvalues of the population covariance matrix. In this note, first, we show that their novel estimator is a
consistent estimator of the population covariance matrix under a high-dimensional asymptotic setup. Moreover, we may
show that the novel estimator is the MLE of the population covariance matrix when $c \in (0, 1)$. The novel estimator is
incorporated to establish the optimal decomposite $T_{T}^{2}-$test for a high-dimensional statistical hypothesis testing
problem and to make the statistical inference for the high-dimensional principal component analysis-related problems
without the sparsity assumption. Some remarks when $p >n $, especially for the high-dimensional low-sample size
categorical data models $p >> n$, are made in the final section.} \\

\vspace{0.3cm}
\small
\parbox{14cm}{\sloppy \, \,
{\bf Keywords:}  High-dimensional covariance matrix, MLEs, the consistent estimator, the decomposite $T_{T}^{2}-$test.}
\small
\parbox{14cm} {\small
}
\end{center}

\vspace{0.3cm}
\def \theequation{1.\arabic{equation}}
\setcounter{equation}{0}
\noindent {\bf 1. Introduction}
\vspace{0.3cm}

\indent The problem in high-dimensional covariance estimation has been one of the most interesting topics in statistics
(Pourahmadi, 2013; Zagidullina, 2021). Stein (1975, 1986) investigated the orthogonally equivariant nonlinear shrinkage
estimator for the population covariance matrix. Stein's estimator has been considered a gold standard, and from which
a large strand of literature on the orthogonally equivariant estimation of covariance matrix was generated (Ledoit and
Wolf, 2012; Rajaratnam and Vincenzi, 2016, and the references therein).

\indent Tsai and Tsai (2024b) also restricted attention to the rotation-equivariant estimators, they showed that the
Stein's estimator can be inadmissible when the dimension $p$ is fixed. Under a high-dimensional asymptotic setup,
namely, both sample size $n$ and the dimension $p$ are sufficiently large with the concentration
$c=\lim_{n \to \infty}{p}/{n}, c \in (0, 1)$, they re-examined the asymptotic optimal property of
estimators proposed by Stein (1975) and Ledoit and Wolf (2018). Moreover, Tsai and Tsai (2024b) looked into the
insight of the Mar$\check{c}$enko-Pastur equation (Silverstein, 1995) to get an explicit equality relationship of the
quantiles of limiting spectral distributions. They used the obtained equality to propose a new kind of orthogonally
equivalent estimator for the population covariance matrix. They showed

\vspace{0.3cm}
\noindent \underline{~~~~~~~~~~~~~~~~~~~~~~~~~~~~~~~~~~~~~~~~~~~~}\\
E-mail: mttsai@stat.sinica.edu.tw
\pagebreak

\noindent  that their novel estimators of the eigenvalues are consistent with the eigenvalues of the population covariance
matrix. When $p/n \to c \in (0, 1)$, they further showed their proposed covariance estimator is the best orthogonally
equivariant estimator for the population covariance matrix under the normalized Stein loss function. In contrast,
both Stein's estimator and the sample covariance matrix can be inadmissible.

\indent The question naturally arises as to whether the consistent estimator of the population covariance matrix exists or
not. In this paper, we further show that the estimator proposed by Tsai and Tsai (2024b) is the consistent estimator
of population covariance matrix $\mbSigma$ when ${p}/{n} \to c \in [0, 1)$. To do that first, we show that the components
for spectral decomposition of the sample covariance matrix are the maximum likelihood estimators (MLEs) of those of the
population covariance matrix when the dimension $p$ is fixed and the sample size $n$ is large (i.e., c=0) in Section 3.
Then, we extend the results of Section 3 to the boundary case, i.e., ${p}/{n} \to c \in (0, 1)$, in Section 4, namely to
show that the novel estimator is not only consistent but also the MLE of the population covariance matrix. Based on the
proposed covariance estimator, the optimal decomposite $T_{T}^{2}-$test for a high-dimensional statistical hypothesis
testing problem is established and it can also be applied to make statistical inferences for the high-dimensional principal
component analysis (PCA) related problems without the sparsity assumption in Section 5. Some remarks when $p > n$,
even for the case $p >> n$, are made in the final section.

\vspace{0.3cm}
\noindent {\bf 2. Preliminary notations}
\def \theequation{2.\arabic{equation}}
\setcounter{equation}{0}
\vspace{0.3cm}

\indent  Let ${\bf X}_{1}, \ldots, {\bf X}_{n}$ be independent $p$-dimensional random vectors with a common
multivariate normal distribution $N_{p}(\bf 0, {\mbSigma})$. A basic problem considered in the literature is the
estimation of the $p \times p$ covariance matrix ${\mbSigma}$, which is unknown and assumed to be non-singular.
It is also assumed that $n \geq p$, as such the sufficient statistic
\begin{align}
{\bf A}=\sum_{i=1}^{n} {\bf X}_{i}{\bf X}^{\top}_{i}
\end{align}
is positive definite with probability one. In the literature, the estimators $\phi({\bf A})$ of ${\mbSigma}$ are the
functions of ${\bf A}$. The sample space ${\mathcal S}$, the parameter space ${\varTheta}$, and the action space
${\mathcal A}$ are taken to be the set $\mathcal P_{p}$ of $p \times p$ symmetric positive definite matrices.
The general linear group $Gl(p)$ acts on the space $\mathcal P_{p}$. Note that ${\bf A}$ has a Wishart distribution
$W({\mbSigma}, n)$, and the maximum likelihood estimator (MLE) of ${\mbSigma}$ is expressed as below
\begin{align}
\Hat{\mbSigma}_{ML}={\bf S}, ~\mbox {where} ~{\bf S}=n^{-1}{\bf A},
\end{align}
which is unbiased (Anderson, 2003).

\indent We consider invariant loss function $L$, i.e., $L$ satisfies the condition that
$L(g\phi({\bf A})g^{\top}, g{\mbSigma}g^{\top})\linebreak =L(\phi({\bf A}), {\mbSigma})$ for all $g \in Gl(p)$.
An estimator $\Hat{\mbSigma}$ is called  $Gl(p)$-equivariant if $\Hat{\mbSigma}({\bf G}{\bf A}{\bf G}^{\top})={\bf G}
\Hat{\mbSigma}({\bf A}){\bf G}^{\top}, \forall {\bf G} \in Gl(p), \forall {\bf A} \in \mathcal P_{p}$. Suppose
that ${\bf G}$ acts on $\mathcal P_{p}$, whereby the orbit through $x \in \mathcal P_{p}$ is the set
${\bf G}x=\{gx|g \in {\bf G}\} \subset \mathcal P_{p}$. The action is called transitive if ${\varTheta}$
is one orbit, i.e., $\forall x, y \in {\varTheta}$ there is some $g \in {\bf G}$ with $gx=y$. It may then
be easy to note the fact that if $L$ is $Gl(p)$-invariant, $\Hat{\mbSigma}$ is $Gl(p)$-equivariant, and
${\bf G}$ acts transitively on $\mathcal P_{p}$, then the risk function is constant on $\mathcal P_{p}$:
$R(\Hat {\mbSigma}, {\mbSigma})=R(\Hat {\mbSigma}, {\bf I}), \forall {\mbSigma} \in \mathcal P_{p}$.

\indent One of the most interesting loss functions was introduced by Stein (1956)
\begin{align}
L(\phi({\bf S}), {\mbSigma})=\mbox{tr}{\mbSigma}^{-1}\phi({\bf S})
      -\mbox{log}\mbox{det}{\mbSigma}^{-1}\phi({\bf S})-p,
\end{align}
where $\mbox{tr}$ and $\mbox{det}$ denote the trace and the determinant of a matrix, respectively. Because
$Gl(p)$ acts transitively on the space $\mathcal P_{p}$, the best $Gl(p)$-equivariant estimator exists.
It can be easily found that the MLE ${\bf S}$ of ${\mbSigma}$ is the best $Gl(p)$-equivariant estimator.
The minimum risk is
\begin{align}
R_{m}(\Hat{\mbSigma}_{ML}, {\mbSigma})=\sum_{i=1}^{p}\{ \mbox{log} n - {\mathcal E}[\mbox{log}{\chi}^{2}_{n-i+1}]\},
\end{align}
where ${\mathcal E}[X]$ denotes the expectation of the random variable of $X$.

%\pagebreak
\vspace{0.3cm}
\noindent {\bf 3. Optimal estimators of ${\mbSigma}$ when it is reparameterized}
\def \theequation{3.\arabic{equation}}
\setcounter{equation}{0}
\vspace{0.3cm}

\indent Since the general linear group $Gl(p)$ is not an amenable group, to study the minimax problem James and Stein (1961)
reparameterized the parameter $\mbSigma \to {\mbTheta}{\mbTheta}^{\top}, {\mbTheta} \in G^{+}_{T}$, where  $G^{+}_{T}$
denotes the group of $p \times p$ lower triangular matrices with positive diagonal elements; the loss function is
also invariant under $G^{+}_{T}$. Using the Cholesky decomposition, we may express that ${\bf A}={\bf T}{\bf T}^{\top}$,
where ${\bf T} \in G^{+}_{T}$. Since $G^{+}_{T}$ acts transitively on the space $\mathcal P_{p}$, the best
$G^{+}_{T}$-equivariant estimator was proposed by James and Stein (1961) as
$\Hat{\mbSigma}_{S}={\bf T}{\bf D}^{-1}_{S}{\bf T}^{\top}$, where ${\bf D}_{S}$ is a positive diagonal matrix with
elements $d_{Sii}=n+p-2i+1, ~i=1, \cdots, p$. The minimum risk for the best $G^{+}_{T}$-equivariant estimator
$\hat{\mbSigma}_{S}$ is
\begin{align}
R_{m}(\Hat {\mbSigma}_{S}, {\mbSigma})=\sum_{i=1}^{p}\{\mbox{log} (n+p-2i+1)
           -{\mathcal E}[\mbox{log}{\chi}^{2}_{n-i+1}]\}.
\end{align}

\indent Because $G^{+}_{T}$ is the solvable group, and hence it is amenable. Thus, Stein's estimator $\Hat{\mbSigma}_{S}$
is minimax.

\vspace{0.3cm}
\noindent {\it 3.1 The Stein phenomenon}
\vspace{0.2cm}

\indent It is easy to see that $ R(\Hat {\mbSigma}_{S}, {\mbSigma}) \leq R(\Hat{\mbSigma}_{ML}, {\mbSigma})$, thus the MLE
${\bf S}$ is inadmissible which people are encouraged to use estimator $\Hat {\mbSigma}_{S}$ instead of ${\bf S}$. This
is the well-known Stein phenomenon for covariance estimation problem, for details see Anderson (2003).

\indent In order to seek the reasons why the Stein phenomenon, which is the MLE ${\bf S}$ of ${\mbSigma}$ inadmissible,
happens. We began to think about the inner meaning of Stein phenomenon. Tsai (2018) extended Stein's method to establish
another minimax estimator. We explain it briefly in the following. Let ${\mbSigma}_{(k)}$ and ${\bf A}_{(k)}$ be
partitioned as
\begin{eqnarray}
  {\mbSigma}_{(k)} = \left[
  \begin{array}{cc}
    {\sigma}_{(k)11} & {\mbSigma}_{(k)12}\\
    {\mbSigma}_{(k)21} & {\mbSigma}_{(k)22}
  \end{array}
  \right]
  ~~\mbox{and}~~
   {\bf A}_{(k)} = \left[
  \begin{array}{cc}
    {a}_{(k)11} & {\bf A}_{(k)12}\\
    {\bf A}_{(k)21} & {\bf A}_{(k)22}
  \end{array}
  \right],
\end{eqnarray}
for all $k=1, \cdots, p$ with ${\mbSigma}_{(1)}={\mbSigma}$ and ${\bf A}_{(1)}={\bf A}$. Define
\begin{align}
{\mbSigma}_{(k+1)}=\mbSigma_{(k)22:1}=\mbSigma_{(k)22}- \mbSigma_{(k)21} \mbSigma_{(k)12}/\sigma_{(k)11}
\end{align}
and
\begin{align}
{\bf A}_{(k+1)}={\bf A}_{(k)22:1}={\bf A}_{(k)22}- {\bf A}_{(k)21}{\bf A}_{(k)12}/a_{(k)11}.
\end{align}
Note that the dimension of ${\mbSigma}_{(k+1)}$ is one less than that of ${\mbSigma}_{(k)}$, which is
a process of successive diagonalization. Let
\begin{align}
{\bf g}_{(k)} = \left[
\begin{array}{cc}
 1  & {\bf 0}    \\
-\mbSigma_{(k)21}\sigma_{(k)11}^{-1} & {\bf I}
\end{array}
\right] ~~\mbox{and}~~
{\bf h}_{(k)} = \left[
\begin{array}{cc}
  1  & {\bf 0}   \\
-{\bf A}_{(k)21}a_{(k)11}^{-1} & {\bf I}
\end{array}
\right], ~k=1, \cdots, p.
\end{align}
We then have:
\begin{align}
\widetilde{\mbSigma}_{(k)} &= {\bf g}_{(k)} {\mbSigma}_{(k)} {\bf g}^{\top}_{(k)}   \\  \nonumber
&= \left[
\begin{array}{cc}
\sigma_{(k)11}& {\bf 0}\\
{\bf 0}& \mbSigma_{(k)22:1}
\end{array}
\right],
\end{align}
and
\begin{align}
\widetilde{\bf A}_{(k)} &= {\bf h}_{k} {\bf A}_{(k)} {\bf h}^{\top}_{(k)}= \left[
\begin{array}{cc}
 a_{(k)11}& {\bf 0}\\
{\bf 0}& {\bf A}_{(k)22:1}
\end{array}
\right], ~~k=1, \cdots, p.
\end{align}
Let
\begin{align}
{\mbSigma}^{*}=\mbox {Diag}({\sigma}_{(1)11}, \cdots, {\sigma}_{(p)11}) ~ \mbox{and}~
{\bf A}^{*}=\mbox {Diag}(a_{(1)11}, \cdots, a_{(p)11}).
\end{align}

\indent Consequently ${\mbSigma}$ and ${\bf A}$ are individually transformed into the diagonal matrices
${\mbSigma}^{*}$ and ${\bf A}^{*}$ so that the one-to-one correspondences of:
${\mbSigma} \leftrightarrow {\mbSigma}^{*}$ and ${\bf A} \leftrightarrow {\bf A}^{*}$ are established to allow
$\phi({\bf A})={\bf D}{\bf A}^{*}$ for Stein loss function, ${\bf D} \in D(p)$, the group of positive diagonal
matrices. By the properties of the Wishart distribution (see Theorem 4.3.4, Theorem 7.3.4, and Theorem 7.3.6 of
Anderson, 2003), it is easy to note that $a_{(i)11}/{\sigma}_{(i)11}, i=1, \cdots, p,$ are independent $\chi^{2}$
random variables with $n-i+1$ degrees of freedom, respectively.  Let ${\bf D}_{0}$ be the diagonal matrix with
elements $d_{0ii}=n-i+1, i=1, \cdots, p$, and we may conclude that ${\bf A}^{*}$ is Wishart distributed with
mean matrix ${\bf D}_{0}{\mbSigma}^{*}$ then. Furthermore, it should be noted that all the $p$ Jacobins of the
transformation of ${\bf A} \to {\bf A}^{*}$ are one, and the Wishart density of ${\bf A}$ is equivalent to the
Wishart density of ${\bf A}^{*}$. Thus the Stein loss function is
\begin{align}
  L(\phi ({\bf A}^{*}), {\mbSigma}^{*})=\mbox{tr}{\mbSigma}^{*-1}{\bf D}{\bf A}^{*}
     -\mbox{log}\mbox{det}{\mbSigma}^{*-1}{\bf D}{\bf A}^{*}-p.
\end{align}
\indent Since ${\bf A}^{*}$ also acts transitively on the space $\mathcal P_{p}$, the best $D(p)$-equivariant
estimator can be expressed as the form of
\begin{align}
  \Hat{\mbSigma}^{*}={\bf D}^{-1}_{0}{\bf A}^{*}.
\end{align}
Thus, the minimum risk for the estimator $\hat{\mbSigma}^{*}_{I}$ is
\begin{align}
  R_{m}(\Hat {\mbSigma}^{*}, {\mbSigma}^{*})=\sum_{i=1}^{p}\{\mbox{log} (n-i+1)
           -{\mathcal E}[\mbox{log}{\chi}^{2}_{n-i+1}]\}.
\end{align}

\indent Since the group $D(p)$ is also solvable, and hence we may conclude that $\Hat {\mbSigma}^{*}$ is a minimax. By (3.1)
and (3.11) it is easy to see that $ R_{m}(\Hat {\mbSigma}^{*}, {\mbSigma}^{*}) \leq  R_{m}(\Hat {\mbSigma}_{S}, {\mbSigma})$,
hence, similar conclusion as the Stein phenomenon we may conclude that Stein's estimator $\Hat {\mbSigma}_{S}$ is
inadmissible, while the estimator $\Hat {\mbSigma}^{*}$ is admissible.

\vspace{0.3cm}
\noindent {\it 3.2 The optimal properties of MLE}
\vspace{0.2cm}

\indent We may note that the MLE ${\bf S}$ of ${\mbSigma}$ is the best $Gl(p)$-equivariant estimator, James and
Stein (1961) used the Cholesky decomposition to parameterize the parameter ${\mbSigma}$ to get the Stein estimator
$\Hat {\mbSigma}_{S}$, which is the best $G_{T}^{+}$-equivariant estimator, Tsai (2018) used the full Iwasawa
decomposition to get the best $D(p)$-equivariant estimator $\Hat {\mbSigma}^{*}$. Note that the inequality
$R_{m}(\Hat {\mbSigma}^{*}, {\mbSigma}^{*}) \leq R_{m}(\Hat {\mbSigma}_{S}, {\mbSigma}) \leq R_{m}
(\Hat {\mbSigma}_{ML}, {\mbSigma})$ holds. Because that $D(p) \subseteq G_{T}^{+} \subseteq Gl(p)$,  we can easily
see that the above inequality holds. The minimum risk of the estimator is larger concerning the larger group and
is smaller concerning the smaller group.

\indent Tsai (2018) showed that the minimum risks of the MLEs under the Cholesky decomposition and the full Iwasawa
decomposition are the same when the geodesic distance loss function on a non-Euclidean space $P_{p}$ is adopted.
Comparing the minimum risks of estimators under different groups does not make much statistical sense. The comparison
of different estimators may make sense when they are compared under the same parameterized decomposition. For the
spectral decomposition, Tsai and Tsai (2024b) claimed that the sample covariance matrix ${\bf S}$ is the best
orthogonally equivariant estimator of spectral decomposition under the Stein loss function. On the other hand, Stein
(1975, 1986) another orthogonally equivariant estimator can be inadmissible under spectral decomposition. These results
contradict the Stein phenomenon that ${\bf S}$ is inadmissible. Hence, the Stein phenomenon seems to be due to the
parameterized decompositions, and it does not seem to have much special statistical meaning.

\indent  Each of the three estimators possesses its optimal properties for their respective parameterized
decomposition. All the three estimators ${\bf S}$, $\Hat {\mbSigma}_{S}$, and $\Hat {\mbSigma}^{*}$ are the best
$Gl(P)$-equivariant, $G_{T}^{+}$-equivariant, and $D(p)$-equivariant estimators, respectively. The sample covariance matrix
${\bf S}$ is not only the best $Gl(p)$-equivariant estimator but also the best $O(p)$-equivariant estimator. They are the
MLEs with respect to $Gl(P)$, $G_{T}^{+}$ and $D(p)$ decompositions, respectively. The optimal property of MLE is essentially
not affected at all. We hope that this paper may impact those statisticians who have been constantly warned not to use
MLE for covariance matrix ever since the Stein phenomenon occurred making them reconsider the employment of the MLE for
covariance matrix.

\indent Note that the Stein loss function is equivalent to the entropy loss function under the multinormal setup. When the
dimension $p$ is fixed and the sample size $n \to \infty$, it has been known in the literature that ${\bf S}$ and
$\Hat {\mbSigma}^{*}$ converge to ${\mbSigma}$ and ${\mbSigma}^{*}$ almost surely $(a.s.)$ when $n \to \infty$,
respectively (Anderson, 2003). We will investigate whether ${\bf S}$ is the MLE of $\mbSigma$ under spectral decomposition
so that the sample components converge to the corresponding population components $a.s.$ as $n \to \infty$, respectively.

\vspace{0.3cm}
\noindent {\it 3.3 The best orthogonally equivariant estimator}
\vspace{0.2cm}

\indent For the application to the statistical inference of principal component analysis, we need the notation of the
so-called spectral decomposition of the population covariance matrix, it can be viewed as another type of reparametrization
of ${\mbSigma}$. Stein (1975, 1986) considered the orthogonally equivariant estimator for the population covariance
matrix, which has been considered a gold standard. Consider the spectral decomposition of the population covariance
matrix, namely ${\mbSigma}={\bf V}{\mbGamma}{\bf V}^{\top}$, where ${\mbGamma}$ is a diagonal matrix with eigenvalues
$\gamma_{i, p}, i=1, \ldots, p$, and ${\bf V}=({\bf v}_1, \ldots, {\bf v}_p)^{\top}$ is the corresponding orthogonal
matrix with ${\bf v}_i$ being the eigenvector associated to the $i$th largest eigenvalue
$\gamma_{i, p}, v_{i1} \geq 0, i=1, \ldots, p$. Similarly, for the sample spectral decomposition, i.e.,
${\bf S}={\bf U}{\bf L}{\bf U}^{\top}$, where ${\bf L}$ is a diagonal matrix with eigenvalues $l_{i, p}$, and
${\bf U}=({\bf u}_1, \ldots, {\bf u}_p)^{\top}$ is the corresponding orthogonal matrix with ${\bf u}_i$ being the
eigenvector corresponding to $l_{i, p}, u_{i1} \geq 0, i=1, \ldots, p$. Write ${\bf L}=\mbox {diag}(l_{1, p},
\ldots, l_{p, p})$ and ${\mbGamma}=\mbox {diag}(\gamma_{1, p}, \ldots, \gamma_{p, p})$. Note that the matrices ${\bf U}$
and ${\bf L}$ are the consistent estimators of ${\bf V}$ and ${\mbGamma}$, respectively when the dimension $p$ is fixed
and the sample size $n$ is large (for the details see Anderson, 2003). Hence, we may conclude that there are two situations
when the dimension $p$ is fixed: (i) When $\mbSigma$ is not reparameterized, the sample covariance matrix ${\bf S}$ is
unbiased and hence it is consistent. (ii) When $\mbSigma$ is reparameterized via spectral decomposition, the the components
${\bf U}$ and ${\bf L}$ are the consistent estimators of ${\bf V}$ and $\mbGamma$, respectively. Then the sample covariance
matrix ${\bf S}$ is still consistent.

\vspace{0.2cm}
\indent {\bf Remark 3.1}. We want to study the consistency property with the help of the optimal properties of MLEs. The
main reason is based on the fact from the general theory of estimation, it is known that the maximum likelihood estimator
is consistent, that is, it tends to the true value with probability one as the sample size becomes large under some
regularity conditions, which are satisfied by the non-degenerated Wishart distribution.

\indent We may note that when $\mbSigma$ is not reparameterized, it is easy to see that the sample covariance matrix
${\bf S}$ is the MLE of ${\mbSigma}$. When the spectral decomposition for ${\mbSigma}$ is adopted, then it is expected
that the sample components ${\bf U}$ and ${\bf L}$ are the MLEs of corresponding population components ${\bf V}$ and
${\mbGamma}$, respectively.

\indent First, when the dimension $p$ is fixed, $n{\bf S}$ is Wishart distributed when $n > p$. Under the spectral
decompositions for ${\mbSigma}$ and ${\bf S}$, we will find the MLEs of ${\bf V}$ and ${\mbGamma}$ in the following.
Note that ${\bf V}, {\bf U} \in \mathcal O(p)$, the set of orthogonal matrices. Let ${\bf H}={\bf V}^{\top}{\bf U}$, then
${\bf H} \in \mathcal O(p)$. Assume that $n \geq p+1$, and then $-\frac{2}{n}$ log-likelihood function of ${\bf S}$ is
\begin{align}
l({\bf S}| {\mbSigma})&=\mbox{tr}{\mbSigma}^{-1} {\bf S}
      -\mbox{log}\mbox{det}{\mbSigma}^{-1}{\bf S}-\frac{2}{n}\mbox{log}c_n({\bf S}) \\ \nonumber
      &=\mbox{tr}{\bf V}{\mbGamma}^{-1}{\bf V}^{\top}{\bf U}{\bf L}{\bf U}^{\top}-\mbox{log}\mbox{det}
      {\bf V}{\mbGamma}^{-1}{\bf V}^{\top}{\bf U}{\bf L}{\bf U}^{\top}-\frac{2}{n}\mbox{log}c_n({\bf L}) \\ \nonumber
      &=\mbox{tr}{\mbGamma}^{-1}{\bf H}{\bf L}{\bf H}^{\top}
      -\mbox{log}\mbox{det}{\mbGamma}^{-1}{\bf L}-\frac{2}{n}\mbox{log}c_n({\bf L}),
\end{align}
where $c_n({\bf S})=\frac{n^{(n-p-1)/2}|{\bf S}|^{-(p+1)/2}}{2^{np/2}\pi^{p(p-1)/4}
\prod_{i=1}^{p}\Gamma[\frac{1}{2}(n-i+1)]}=c_n({\bf L})$ which is independent of ${\mbSigma}$ (i.e., ${\bf V},
{\mbGamma}$). The equation (3.12) is essentially equivalent to the Stein loss function.

\vspace{0.3cm}
\indent {\bf Theorem} (von Neumann, 1937). {\it For ${\bf H}$ orthogonal and ${\bf D}_{\mbgamma}$ and
${\bf D}_{\bf l}$ diagonal} ($\gamma_1 \geq \ldots \geq \gamma_p > 0, l_1 > \ldots > l_p > 0$)
\begin{align}
 \mbox{min}_{{\bf H}\in \mathcal O(p)}~\mbox{tr}{\bf D}^{-1}_{\mbgamma}{\bf H}{\bf D}_{\bf l}{\bf H}^{\top}
 =\mbox{tr}{\bf D}^{-1}_{\mbgamma}{\bf D}_{\bf l},
\end{align}
{\it and a minimizing value of ${\bf H}$ is $\Hat{\bf H}={\bf I}$}. For the details of proofs, see Theorem A.4.7
and Lemma A.4.6 of Anderson (2003).
\vspace{0.2cm}

\indent By the result of von Neumann Theorem, we then have the MLE of $V$ is that $\Hat{\bf V}={\bf U}$, and hence
\begin{align}
&\mbox{min}_{{\bf V} \in \mathcal O(p)}\mbox{tr}{\bf V}{\mbGamma}^{-1}{\bf V}^{\top}{\bf U}
  {\bf L}{\bf U}^{\top} \\ \nonumber
%&=\mbox{min}_{{\bf H} \in \mathcal Q(p)}\mbox{tr}{\mbGamma}^{-1}{\bf H}\Hat\mbPsi({\bf L}){\bf H}^{\top}  \\ \nonumber
 & =\mbox{tr}{\mbGamma}^{-1}{\bf L}.
\end{align}
Thus, we may have
\begin{align}
\mbox{min}_{{\bf V} \in \mathcal O(p)}l({\bf S}| {\mbSigma})=\mbox{tr}{\mbGamma}^{-1}{\bf L}
      -\mbox{log}\mbox{det}{\mbGamma}^{-1}{\bf L}-\frac{2}{n}\mbox{log}c_n({\bf L}).
\end{align}

\indent After some calculations, the function $\mbox{min}_{{\bf V}\in {\mathcal O}(p)}l({\bf S}|{\mbSigma})$ in (3.15)
is further minimized with respect to ${\mbGamma}$ at $\Hat{\mbGamma}={\bf L}$. As such, when $p$ is fixed, we have
that ${\bf U}$ is the MLE of ${\bf V}$ and $l_{i, p}$ is the MLE of $\gamma_{i, p}, i=1, \ldots, p$. Thus, when
$p$ is fixed and the sample size $n$ is large, according to the property of MLE we have ${\bf U} \to {\bf V}$ and
${\bf L} \to {\mbGamma}$ almost surely (a.s.), and thus ${\bf U}$ and ${\bf L}$ are the consistent estimators of
${\bf V}$ and ${\mbGamma}$, respectively. Hence, ${\bf S}={\bf U}{\bf L}{\bf U}^{\top} \to {\bf V}{\mbGamma}{\bf V}^{\top}
={\mbSigma}$ a.s. Therefore, in terms of the spectral decompositions, when the dimension $p$ is fixed the sample covariance
matrix ${\bf S}$ is the consistent estimator of the population covariance matrix ${\mbSigma}$. From the above arguments,
when the dimension $p$ is fixed we may note that the MLEs play an important role in being optimal whether it is
reparameterized or not. When it is not reparameterized, the MLE ${\bf S}$ of ${\mbSigma}$ is unbiased and consistent,
while it is reparameterized, the MLEs of component parameters for spectral decomposition are consistent.

\indent However, the situation may be different because the sample covariance matrix ${\bf S}$ is no longer to be the
MLE of the population covariance matrix ${\mbSigma}$ anymore when the dimension $p$ is large such that $c \in (0, 1)$.
Hence, the question naturally arises as to whether the consistent estimator of ${\mbSigma}$ exists or not under the
large dimensional asymptotics setup. Under the spectral decomposition, Tsai and Tsai (2024b) proved the consistency
for their proposed estimators of population eigenvalues with the help of random matrix theory. Some notations of it
are presented below.

\vspace{0.3cm}
\noindent {\bf 4. High-dimensional case}
\def \theequation{4.\arabic{equation}}
\setcounter{equation}{0}
\vspace{0.3cm}

\indent For a large $(n, p)$ set up, the large dimensional asymptotics framework is setted up when $(n, p) \to \infty$
such that $c=\lim_{n \to \infty}p/n$ is fixed, $0 \leq c < 1$. In this section, we extend the class of orthogonally
equivariant estimators to the realm of large dimensional asymptotics with the concentration $c \in (0,1)$.

\vspace{0.3cm}
\noindent {\it 4.1 The Mar$\check{c}$enko-Pastur equation}
\vspace{0.2cm}

\indent The same as Ledoit and P$\acute{e}$ch$\acute{e}$ (2011), we make the following assumptions:

\indent A1. Note that ${\bf x}_{i}={\mbSigma}^{1/2}{\bf z}_{i}, i=1,  \ldots, n$, where ${\bf z}_{i}$ are independent and
identically distributed with mean $ {\bf 0}$ and covariance matrix ${\bf I}$. Assume that the 12th absolute central
moment of each variable $z_{ij}$ bounded by a constant.

\indent A2. The population covariance matrix $ {\mbSigma}$ is nonrandom positive definite.
$\mbox{lim inf}_{p \to \infty} \gamma_{p, p} > 0$ and $\mbox{lim sup}_{p \to \infty} \gamma_{1, p} < \infty$.

\indent A3. For large $(n,p)$ set-up, the large dimensional asymptotics framework is setted up when $(n,p) \to \infty$
such that $c = p/n$ is fixed $0 \leq c < 1$ in this paper.

\indent A4. Let $0 < \gamma_{p, p} < \cdots < \gamma_{1, p} $. The emperical spectral distribution of $\mbSigma$
defined by $H_{n}(\gamma)=\frac{1}{p}\sum_{i=1}^{p}1_{[\gamma_{i, p}, \infty)}(\gamma)$, converges as $p \to \infty$ to a
probability distribution function $H(\gamma)$ at every point of continuity of $H$. The support of $H$, $\mbox{Supp}(H)$,
is included in a compact set $[h_1, h_2]$ with $0 < h_1 \leq h_2 < \infty$.

\indent Let $F_{n}(\lambda)=\frac{1}{p}\sum_{i=1}^{p}1_{[l_{i, p}, \infty)}(\lambda)$ be the sample spectral
distribution and $F$ be its limiting. It is proved that $F_{n}$ converges to $F$ a.s. as $n \to \infty$
(Mar$\check{c}$enko-Pastur, 1967).

\indent The Stieltjes transform of distribution function $F$ is defined by
\begin{align}
m_F(z)=\int_{-\infty}^{\infty}\frac{1}{l-z}dF(l), \forall z \in C^{+},
\end{align}
where $C^{+}$ is the half-plane of complex numbers with a strictly positive imaginary part.
Let
\begin{align}
m_{F_{n}}(z)=p^{-1}\mbox{tr}[({\bf S}-z{\bf I})^{-1}],
\end{align}
then from the results of random matrix theory $F_n(z)$ converges to $F(z)$ if and only if $m_{F_{n}}(z)$
converges to $m_F(z)$. Subsequently, the well-known Mar$\check{c}$enko-Pastur equation (Silverstein, 1995) can be
expressed in the following form
\begin{align}
m_F(z)=\int_{-\infty}^{\infty}\frac{1}{\gamma[1-c-czm_F(z)]-z}dH(\gamma), \forall z \in C^{+},
\end{align}
where $H$ denotes the limiting behavior of the population spectral distribution. Upon the Mar$\check{c}$enko-Pastur
equation, meaningful information of the population spectral distribution can be retrieved under the large
dimensional asymptotics framework. Choi and Silverstein (1995) further showed that
\begin{align}
\lim_{z\in C^{+}\to \l}m_F(z)=\check{m}_F(l)
\end{align}
exists for any $l \in R / \{0\}$.

\indent Using the Sokhotski-Plemelj formula, the term $\check{m}_F(l)$ can be seperated into the real part which
becomes a principal value integral (the so-called Hilbert transform), while the imaginary part becomes $\pi$ times
the limiting sample spectral density function $f(l)$. Namely,
\begin{align}
\check{m}_F(l)=\mbox{Re}[\check{m}_F(l)]+i\pi f(l),
\end{align}
where the Hilbert transform denotes
\begin{align}
\mbox{Re}[\check{m}_F(x)]=Pr\int \frac{dF(t)}{t-x}.
\end{align}

\indent For some special cases, $\check{m}_F(x)$ can be expressed explicitly. For example,
let $\lambda_{+}=(1+\sqrt{c})^{2}$ and $\lambda_{-}=(1-\sqrt{c})^{2}$. When $\mbSigma={\bf I}$,
then the Mar$\check{c}$enko-Pastur density function is of the form
\begin{align}
f_{MP}(x)=\frac{\sqrt{(x-\lambda_{-})(\lambda_{+}-x)}}{2\pi cx}, ~~x \in (\lambda_{-}, \lambda_{+}).
\end{align}
By the resolvent method, we then have
\begin{align}
\check{m}_F(x)=\frac{1-c-x+i\sqrt{(x-\lambda_{-})(\lambda_{+}-x)}}{2cx},
\end{align}
where the real part is the Cauchy principal value, i.e.,
\begin{align}
\mbox{Re}[\check{m}_F(x)]&=Pr\int f_{MP}(t)\frac{dt}{t-x}  \\ \nonumber
&=\frac{1-c-x}{2cx}.
\end{align}
Generally, ${\mbSigma}$ is unknown, and the form of $\mbox{Re}[\check{m}_F(x)]$ will not be explicit,

\indent Stein (1975) used the naive empirical counter part $\check{m}_{F_{n}}(l_{i, p})(=\frac{1}{p}\sum_{j \ne i}
\frac{1}{l_{j, p}-l_{i, p}})$ to estimate the Hilbert transformation $\mbox{Re}[\check{m}_F(l_{i})]$, where $l_{i}$
denotes the $(1-\alpha)$-quantile of limiting sample spectral distribution $F$ so that $[p(1-\alpha)]=i, i=1, \ldots, p$,
with $[x]$ denoting the largest integer of $x$. Since that $F_n(z)$ converges to $F(z)$ a.s., as such $m_{F_{n}}(z)$
converges to $m_F(z)$ a.s. Thus we have that $l_{i, p}$ converges to $l_{i}~a.s., ~i=1, \ldots, p$. And then the
estimator $\check{m}_{F_{n}}(l_{i, p})$ proposed by Stein is a consistent estimator of $\mbox{Re}[\check{m}_F(l_{i})]$.

\vspace{0.3cm}
\noindent {\it 4.2 The consistent estimators of population eigenvalues}
\vspace{0.2cm}

\indent  The Mar$\check{c}$enko-Pastur equation in (4.3) shows the implicit relationship between $F$ and $H$, Tsai and
Tsai (2024b) further established the following explicit equality relationship
\begin{align}
\gamma_{i}= \frac{l_{i}}{1-c-cl_{i}\mbox{Re}[\check{m}_F(l_{i})]},  i=1, \ldots, p,
\end{align}
where $\gamma_{i}$ and $l_{i}$ denote the $(1-\alpha)$-quantiles of limiting population and sample spectral distributions
$H$ and $F$, respectively, so that $[p(1-\alpha)]=i, i=1, \ldots, p$, with $[x]$ denoting the largest integer of $x$.
Let $\mbox{Supp}(F)$ be the support of $F$. Via Theorem 4.1 of Choi and Silverstein (1995), Ledoit and
P$\acute{e}$ch$\acute{e}$ (2011) pointed out that if $l_{i} \notin \mbox{Supp}(F)$, then
${l_{i}}/{1-c-cl_{i}\mbox{Re}[\check{m}_F(l_{i})]} \notin \mbox{Supp}(H),$ for $l_{i} \in R/\{0\}, i=1, \ldots, p.$

\indent Write $\gamma_{i}=\psi_{i}({\bf L}), i=1, \ldots, p$, Tsai and Tsai (2024b) proposed a new kind of orthogonally
equivariant estimator $\Hat {\mbSigma}_{T}$ of $\mbSigma$, which is of the form
\begin{align}
\Hat {\mbSigma}_{T}={\bf U}\Hat\mbPsi({\bf L}){\bf U}^{\top}, &~ \mbox {where}~
\Hat\mbPsi({\bf L})=\mbox {diag}(\Hat\psi_{1}({\bf L}), \cdots, \Hat\psi_{p}({\bf L})) ~\mbox {with}  \\ \nonumber
\Hat\psi_{i}({\bf L})&=\frac{nl_{i, p}}{n-p+1-pl_{i, p}{\check{m}}_{F_n}(l_{i, p})}  \\ \nonumber
         &=nl_{i, p}(n-p+1-l_{i, p}\sum_{j \ne i}\frac{1}{l_{j, p}-l_{i, p}})^{-1}, ~ i=1, \ldots, p.
\end{align}

\indent When $c=0$ as discussed in Section 3, we have ${\gamma}_{i, p} \to {\gamma}_{i}$ and
${\gamma}_{i}=l_{i}, i=1, \ldots, p$. However, when $p/n \to c \in (0, 1)$, $\gamma_{i, p} \to {\gamma}_{i}$ and
$\gamma_{i}$ is no longer to be $l_{i}$ any more. But via the equation (4.10), it should be of the form
$\frac{l_{i}}{1-c-cl_{i}[\mbox{Re}\check{m}_F(l_{i})]}, i=1, \ldots, p.$ Note that by the assumption A4 that $H_{n}$
converges to $H$ when $c \in (0, 1)$, thus ${\gamma}_{i, p}$ converges to ${\gamma}_{i}, i=1, \ldots, p$, namely
${\mbGamma}$ converges to $\mbPsi({\bf L})$ defined in Proposition 4.1. Hence, to estimate ${\gamma}_{i, p}$ is the
same as to estimate ${\gamma}_{i}$ under the large dimensional asymptotics setup, $i=1, \ldots, p$. Under some regularity
conditions, Tsai and Tsai (2024b) claimed that their proposed estimators of the population eigenvalues are consistent.
We summarize the results in the following.

\vspace{0.3cm}
\indent {\bf Proposition 4.1}. {\it Let $\mbPsi({\bf L})=\mbox {diag}(\gamma_{1}, \ldots, \gamma_{p})$ and
 $\Hat\mbPsi({\bf L})=\mbox {diag}(\Hat\psi_{1}({\bf L}), \ldots, \Hat\psi_{p}({\bf L}))$ be difined in (4.10) and
(4.11), respectively. Under the assumptions of Theorem 4.1 of Tsai and Tsai (2024b), then $\Hat\mbPsi({\bf L})$ is the
consistent estimator of $\mbPsi({\bf L})$, namely $\Hat\mbPsi({\bf L})$ is the consistent estimator of $\mbGamma$, when
$p/n \to c \in (0, 1)$.}
\vspace{0.2cm}

\indent {\bf Remark 4.1}. Under some regularity conditions, the new explicit equality relationship between the quantiles
of limiting sample spectral distribution $F$ and limiting population spectral distribution $H$ is established in (4.10),
as such, the consistent estimators of the population eigenvalues can then be easily found. This result makes up for the
deficiency of both Stein's estimator and Ledoit and Wolf's estimator, which are not consistent estimators of population
eigenvalues. Random matrix theory did the essential help of our finding. However, it remains unsolved whether the estimator
$\Hat \mbSigma_{T}$ is consistent for $\mbSigma$, namely, whether the sample component estimators are consistent for the
corresponding population components. To overcome this difficulty, we will adopt the MLE approach to investigate it.

\vspace{0.3cm}
\noindent {\it 4.3 The consistent estimator of the population covariance matrix}
\vspace{0.2cm}

\indent When the dimension $p$ is large, the sample covariance matrix ${\bf S}$ is no longer to be the MLE of ${\mbSigma}$
any more. It is difficult to directly find out the functional form of ${\bf S}$ so that it is the MLE of ${\mbSigma}$, as
such we may take the detour of the reparameterization of ${\mbSigma}$ via spectral decomposition to overcome the difficulty.
The main goal next is to see whether the orthonormal matrix ${\bf U}$ is the MLE of ${\bf V}$ or not.
The result of ${\mathcal E}{\bf U}={\bf V}$ implies that the limiting distribution of ${\bf U}$ on ${\mathcal O}(p)$
is entirely concentrated at ${\bf V}$, the unbiasedness is not a useful optimal property, the role of which would be
replaced by the property of equivariance. Ledoit and P$\acute{e}$ch$\acute{e}$ (2011) pointed out that the projection
of the sample eigenvector onto the population eigenvector $p|{\bf u}^{\top}_{i}{\bf v}_{j}|^{2}$ will wipe out the
non-rotation equivariant behavior and the average of the quantities of $p|{\bf u}^{\top}_{i}{\bf v}_{j}|^{2}$ over the
sample eigenvectors associated with the sample eigenvalues, quantities how the eigenvectors of the sample covariance
matrix deviate from those of the population covariance matrix under the large-dimensional asymptotics. This is one of
the main reasons we prefer to restrict it to the class of rotation-equivariant estimators. Tsai and Tsai (2024b)
established the best orthogonally equivariant estimator $\Hat {\mbSigma}_{T}$ for ${\mbSigma}$. We continue to study
whether the proposed estimator $\Hat {\mbSigma}_{T}$ is the consistent estimator of the population covariance matrix
${\mbSigma}$ or not when $p/n \to c\in (0,1)$. By Proposition 4.1, it only needs to see whether ${\bf U}$ is the consistent
estimator of ${\bf V}$.

\indent The orthogonal matrix ${\bf U}$ may not generally be a consistent estimator of ${\bf V}$ when the dimension $p$ is
large (see Bai et al., 2007, and and references therein). Hence, we may work it under the restricted model, namely, under
the Wishart distribution setup when $p/n \to c \in (0, 1)$.

\indent When ${\mbSigma}$ is reparameterized via the spectral decomposition, we want to study the consistency property
of component parameters when the dimension $p$ is large. Under the multivariate normal setup, when the dimension $p$ is
fixed, $n{\bf S}$ is Wishart distributed with the mean matrix $\mbSigma$. However, when $\lim_{n \to \infty} p/n=c
\in (0, 1)$, $n{\bf S}$ is neither to be the Wishart distributed with mean matrix ${\mbSigma}$, nor ${\bf S}$ is the MLE
of ${\mbSigma}$. Instead, we may notice that $n{\Hat \mbSigma}_{T}$ is Wishart-type distributed with mean matrix
${\mbSigma}$ when $p/n \to c\in (0,1)$. It is easy to note that the $-\frac{2}{n}$ log-likelihood function
$l(\Hat {\mbSigma}_{T}|{\mbSigma})$ of ${\Hat\mbSigma}_{T}$ is similar to (3.12), with $\Hat{\mbSigma}_{T}$ in (4.11)
replacing ${\bf S}$ in (3.12) (i.e., with $\Hat{\mbPsi({\bf L})}$ replacing ${\bf L}$), which still satisfies the regularity
conditions, and it does not degenerate. Based on $l(\Hat {\mbSigma}_{T}|{\mbSigma})$, our goal is to show  that
$\Hat{\mbSigma}_{T}$ is the MLE of ${\mbSigma}$ when $\lim_{n \to \infty} p/n=c\in (0, 1)$.

\indent First, we want to show that ${\bf U}$ is the MLE of ${\bf V}$ when $p/n \to c\in (0, 1)$, namely to extend the
von Neumann Theorem to the case $p/n \to c\in (0, 1)$. Note that ${\bf H}{\bf H}^{\top}={\bf I}$, thus
$d{\bf H}{\bf H}^{\top}+{\bf H}d{\bf H}^{\top}={\bf 0}$. Moreover,
\begin{align}
d\mbox{tr}{\mbGamma}^{-1}{\bf H}\Hat\mbPsi({\bf L}){\bf H}^{\top}
&=\mbox{tr}{\mbGamma}^{-1}d{\bf H}\Hat\mbPsi({\bf L}){\bf H}^{\top}+\mbox{tr}{\mbGamma}^{-1}{\bf H}
\Hat\mbPsi({\bf L})d{\bf H}^{\top} \\ \nonumber
&=\mbox{tr}{\mbGamma}^{-1}d{\bf H}\Hat\mbPsi({\bf L}){\bf H}^{\top}-\mbox{tr}{\mbGamma}^{-1}{\bf H}
\Hat\mbPsi({\bf L}){\bf H}^{\top}d{\bf H}{\bf H}^{\top},
\end{align}
then the derivative becomes $\frac{d\mbox{tr}{\mbGamma}^{-1}{\bf H}\Hat\mbPsi({\bf L}){\bf H}^{\top}}{d{\bf H}}
={\mbGamma}^{-1}\Hat\mbPsi({\bf L}){\bf H}^{\top}-{\mbGamma}^{-1}{\bf H}\Hat\mbPsi({\bf L}){\bf H}^{\top}{\bf H}^{\top}$.
Thus, $\frac{d\mbox{tr}{\mbGamma}^{-1}{\bf H}\Hat\mbPsi({\bf L}){\bf H}^{\top}}{d{\bf H}}={\bf 0}$ implies
that ${\bf H}\Hat\mbPsi({\bf L}){\bf H}^{\top}=\Hat\mbPsi({\bf L})$. Similar arguments as above, we can also show that
${\bf H}{\mbGamma}^{-1}{\bf H}^{\top}={\mbGamma}^{-1}$. Hence we may have that
${\mbox min}_{{\bf H}\in {\mathcal O(p)}}{\mbox tr}{\mbGamma}^{-1}{\bf H}\Hat\mbPsi({\bf L}){\bf H}^{\top}
={\mbox tr}{\mbGamma}^{-1}\Hat\mbPsi({\bf L})$, namely the minimum of $\mbox{tr}{\mbGamma}^{-1}{\bf H}\Hat\mbPsi({\bf L})
{\bf H}^{\top}$, with respective to ${\bf H}\in {\mathcal O(p)}$, occurs at $\Hat{\bf H}={\bf I}$. Thus, the von Neumann
Theorem still holds for the boundary case, i.e., when $p/n \to c \in (0, 1)$. As such, in terms of the spectral
decompositions, we may have that ${\bf U}$ is also the MLE of ${\bf V}$ when $p/n \to c \in (0, 1)$. Hence, by the
property of MLE we may summarize it as in the following.

\vspace{0.3cm}
\indent {\bf Theorem 4.1}. {\it Let ${\bf X}_{1}, \ldots, {\bf X}_{n}$ be independent $p$-dimensional random vectors
with a common multivariate normal distribution $N_{p}(\bf 0, {\mbSigma})$. Consider the spectral decompositions
$\mbSigma={\bf V}{\mbGamma}{\bf V}^{\top}$ and ${\bf S}={\bf U}{\bf L}{\bf U}^{\top}$, and let $\Hat {\mbSigma}_{T}$
be defined as in (4.11). Under the assumptions of Proposition 4.1. When $\lim_{n \to \infty} p/n=c \in (0, 1)$, then we
have that ${\bf U}$ is the MLE of ${\bf V}$, hence, it is the consistent estimator of ${\bf V}$.}
\vspace{0.2cm}

\indent By Proposition 4.1 and Theorem 4.1, we then may conclude that the proposed novel estimator $\Hat {\mbSigma}_{T}$
is consistent for the population covariance matrix ${\mbSigma}$ when the dimension $p$ is large. Next, we continue to
investigate whether the proposed estimator $\Hat {\mbSigma}_{T}$ is the MLE of $\mbSigma$ or not. Note that
\begin{align}
\mbox{min}_{{\bf V} \in \mathcal O(p)}l(\Hat{\mbSigma}_{T}| {\mbSigma})=\mbox{tr}{\mbGamma}^{-1}\Hat\mbPsi({\bf L})
      -\mbox{log}\mbox{det}{\mbGamma}^{-1}\Hat\mbPsi({\bf L})-\frac{2}{n}\mbox{log}c_n(\Hat\mbPsi({\bf L})).
\end{align}
After some calculations the function $\mbox{min}_{{\bf V}\in {\mathcal O}(p)}l(\Hat{\mbSigma}_{T}|{\mbSigma})$
in (4.13) is further minimized with respect to ${\mbGamma}$ at $\Hat{\mbGamma}=\Hat\mbPsi({\bf L})$. As such, when
$p/n \to c\in (0,1)$, we may conclude that $\Hat\mbPsi({\bf L})$ is the MLE of $\mbSigma$, however, the sample covariance
matrix ${\bf S}$ is not. According to the property of MLE we have that ${\bf U}$, $\Hat\mbPsi({\bf L})$ and
$\Hat \mbSigma_{T}$ are the consistent estimators of ${\bf V}$, $\mbGamma$, and ${\mbSigma}$, respectively. Therefore,
we have the following.

\vspace{0.3cm}
\indent {\bf Theorem 4.2}. {\it Under the assumptions of Theorem 4.1.  For the boundary case, i.e., when
$ p/n \to c \in (0, 1)$), then $\Hat {\mbSigma}_{T}$ is the MLE of ${\mbSigma}$. Hence, it is consistent.}
\vspace{0.2cm}

\indent {\bf Remark 4.2}. We may conclude three-fold in the following: (i) The sample covariance matrix ${\bf S}$ is the
MLE of population covariance matrix $\mbSigma$ when the dimension $p$ is fixed. (ii) The estimator $\Hat {\mbSigma}_{T}$
is the MLE of ${\mbSigma}$ when the dimension $p$ is large such that $\lim_{n \to \infty} p/n=c \in (0, 1)$. (iii) It is
easy to see that $\Hat {\mbSigma}_{T}$ reduces to the sample covariance matrix ${\bf S}$ when the dimension $p$ is fixed
and the sample size $n$ is large (i.e., $c=0$). Those are insightful parallels. Hence, for simplicity, we may integrate
the above results into a unified one: when $p$ is fixed or $\lim_{n\to \infty}p/n=c \in (0, 1)$ (i.e., $c \in [0, 1))$,
$n{\Hat \mbSigma}_{T}$ is Wishart distributed with mean matrix ${\mbSigma}$ and ${\Hat \mbSigma}_{T}$ is the MLE of
${\mbSigma}$. Thus, ${\Hat \mbSigma_{T}}$ is the consistent estimator of ${\mbSigma}$, hence, ${\Hat \mbSigma}_{T}$
converges to ${\mbSigma}~a.s.$ as ${n \to \infty}.$ Therefore, we may use ${\Hat \mbSigma}_{T}$ to replace ${\bf S}$ for
making statistical inferences when the dimension $p$ is fixed or $c \in (0, 1))$.

\vspace{0.2cm}
\indent {\bf Remark 4.3}. Tsai and Tsai (2024b) used the fundamental statistical concept to get the quantile equality
relationship of limiting sample and population spectral distributions so that the consistent problems between the sample
eigenvalues and the population eigenvalues can be easily handled. Then use the notion of the likelihood function to get
things done. As long as the density function does not degenerate, the statistical inference can be performed similarly
to the traditional one. The key point in having this conclusion is to find a consistent estimator of the population
covariance matrix. Namely, it is directly to find out the MLE of ${\mbSigma}$ when $n > p$.

\indent When the dimension $p$ is fixed, we have that $l_{i, p}$ is the MLE of $\gamma_{i, p}, i=1, \ldots, p$. However,
 it is not true that $l_{i, p}$ is the MLE of $\gamma_{i, p}, i=1, \ldots, p$, when $p/n \to c \in (0, 1)$.

\indent Johnstone and Paul (2018) provided a detailed discussion on sample eigenvalue bias and eigenvector
inconsistency under the spiked covariance model and for high-dimension PCA-related phenomena.

\vspace{0.2cm}
\indent {\bf Remark 4.4}. We may note that $n\Hat {\mbSigma}_{T}$ is Wishart distributed with the mean matrix ${\mbSigma}$,
when the dimension $p$ is fixed, then $\Hat {\mbSigma}_{T}$ reduces to the sample covariance matrix ${\bf S}$. With the
proposed novel estimator $\Hat {\mbSigma}_{T}$ replacing the sample covariance matrix ${\bf S}$ for statistical inference,
the case for the traditional fixed dimension $p$ and the case for the nowadays high-dimensional can be integrated into one.
From the above arguments, we may suggest using the proposed consistent estimator $\Hat{\mbSigma}_{T}$ to replace the sample
covariance matrix ${\bf S}$ to make the multivariate statistical inference including the PCA-related problems for the cases
as long as (i) when $p$ is fixed and $n$ is large , i.e., $c=0$, and (ii) when $p/n \to c \in (0, 1)$.

\vspace{0.2cm}
\indent We provide an outline for the likelihood ratio test (LRT) of the hypothesis testing problem in the next section.

\vspace{0.3cm}
\noindent {\bf 5. The decomposite $T_{T}^{2}-$test when the dimension $p$ is large}
\def \theequation{5.\arabic{equation}}
\setcounter{equation}{0}
\vspace{0.3cm}

\indent Let ${\bf X}_{i}, i=1, \ldots, n,$ be $n$ $i.i.d$ random vector having a $p$-dimensional multinormal distribution
with mean vector ${\mbmu}$ and unknown positive definite covariance matrix ${\mbSigma}$. Consider the hypothesis testing
problem
\begin{align}
H_{0}:{\mbmu}={\bf 0}~\mbox{versus}~H_{1}:{\mbmu}\ne {\bf 0}
\end{align}
when both the dimension $p$ and the sample size $n$ are large. Let
\begin{align}
{\mXbar}=\frac{1}{n}\sum_{i=1}^{n}{\bf X}_{i} ~\mbox{and}~{\bf S}=\frac{1}{n-1}\sum_{i=1}^{n}
({\bf X}_{i}-{\mXbar})({\bf X}_{i}-{\mXbar})^{\top},
\end{align}
then the well-known Hotelling's $T^{2}$-test statistic in the literature is denoted as
\begin{align}
T^{2}=n{\mXbar}^{\top}{\bf S}^{-1}{\mXbar}.
\end{align}
When the dimension $p$ is fixed, Hotelling's $T^{2}$-test is optimal for the problem (5.1).

\indent However, when the dimension $p$ is large, the performance of Hotelling's $T^{2}$-test is not optimal due to
the fact that the sample covariance matrix ${\bf S}$ is no longer to be the consistent estimator of ${\mbSigma}$.
To overcome the difficulty, we may adopt the novel estimator $\Hat{\mbSigma}_{T}$ to replace the sample covariance
matrix ${\bf S}$, and then consider the following decomposite $T_{T}^{2}$-test statistic
\begin{align}
T^{2}_{T}=n{\mXbar}^{\top}\Hat{\mbSigma}_{T}^{-1}{\mXbar},
\end{align}
where $\Hat{\mbSigma}_{T}$ is defined in (4.11) with ${\bf S}$ in (5.2) replacing the one in (2.2). It is easy to note that
$T_{T}^{2}$-test is the LRT statistic for the problem (5.1).

\indent Since the power of any reasonable test goes to one as $n \to \infty$, to avoid the difficulty Le Cam's contiguity
concept was adopted to study the asymptotically local distribution when the dimension $p$ is fixed.  Note that the
traditional local alternatives do not depend on the dimension $p$. When the dimension $p$ is large, Chia-Hsuan in her Ph.D.
thesis (Tsai and Tsai, 2024a) incorporated the dimension $p$ into the consideration to study asymptotical distribution
under the local alternatives
\begin{align}
H_{0}:{\mbmu}={\bf 0}~\mbox{versus}~H_{1n}:{\mbmu}=n^{-1/2}p^{1/4}{\mbdelta},
\end{align}
where ${\mbdelta}$ is a fixed $p$-dimensional vector, which means to assume that ${\mbdelta}^{\top}{\mbSigma}^{-1}
{\mbdelta} < \infty$ when $p$ is large. Comparing to the traditional one, the local alternatives also depend on the
dimension $p$ with a little bit of change of converge rate. Let
\begin{align}
T^{2}_{0}=n{\mXbar}^{\top}{\mbSigma}^{-1}{\mXbar}.
\end{align}
Then similar arguments as that of Tsai and Tsai (2024a), when $p/n \to c \in [0, 1)$ it can be shown that $T^{2}_{T}$ does
not converge to $T^{2}_{0}$ in probability, however, it is true that $T^{2}_{T}$ converges to $T^{2}_{0}$ in distribution
locally. Note that in the traditional case of the fixed dimension $p$, the proposed decomposite $T^{2}_{T}$-test reduces
to Hotellin's $T^{2}$-test, and Hotellin's $T^{2}$-test statistic converges to $T^{2}_{0}$ in probability, which implies
the convergence in distribution. It is not hard to see that $T_{T}^{2}-$test statistic asymptotically locally (under $H_{1n}$ with the rate
$n^{-1/2}p^{1/4}$) reduces to non-central chi-square ${\chi}^{2}_{p}({\mbdelta}^{\top}{\mbSigma}^{-1}{\mbdelta})$
distributed. This asymptotically local power function is still the monotone function of non-centrality ${\mbdelta}^{\top}
{\mbSigma}^{-1}{\mbdelta}$. Hence, when $p/n \to c \in [0, 1)$ it is easy to see that the proposed decomposite
$T^{2}_{T}$-test is optimal for the problem (5.1), for details, see Tsai and Tsai (2024a).

\vspace{0.2cm}
\indent {\bf Remark 5.1}. The high-dimensional PCA problem has been mainly studied under the spiked covariance models,
there is a need to make the sparsity assumption on the population eigenvectors for the consistent problem (Johnstone
and Lu, 2009, and the references therein). On the other hand, with the proposed novel estimator $\Hat {\mbSigma}_{T}$
replacing the sample covariance matrix ${\bf S}$ for statistical inference, when $p/n \to c \in [0, 1)$, the results
of Theorem 4.2 can be applied to make multivariate statistical inferences and the PCA-related problems without the
sparsity assumption. When $c \in [0, 1)$, our approach unifies the traditional case ($c=0$) and modern high-dimensional
case $(c \in (0, 1))$ for the multivariate statistical methods and high-dimensional PCA-related problems. The proposed
novel estimator is incorporated to establish the optimal decomposite $T_{T}^{2}-$test for a high-dimensional statistical
hypothesis testing problem and can be directly applied to high-dimensional PCA-related problems without the sparsity
assumption.

\indent Some remarks when $p >n $, especially for the high-dimensional low-sample size categorical data models $p >> n$,
are made in the final section.

\vspace{0.3cm}
\noindent {\bf 6. General remarks for the situation when $p > n$ }
\def \theequation{6.\arabic{equation}}
\setcounter{equation}{0}
\vspace{0.3cm}

\noindent {\it 6.1. When $p > n$, both $n$ and $p$ are fixed}
\vspace{0.2cm}

\indent Under the multivariate normal setup, we may note that when $p > n$, $p$ and $n$ are fixed, then the density
function of ${\bf S}$ becomes the singular Wishart distribution (Uhlig, 1994), which it degenarates. In this situation,
assume that rank(${\bf S})=n$ we then may have the notations: ${\bf L}$ is a $n \times n$ diagonal matrix, and the
reparameterization ${\bf S}={\bf U_{1}}{\bf L}\bf U_{1}^{\top}$, where ${\bf U_{1}}\in {\mathcal V_{n, p}}$, the
$np-n(n+1)/2$-dimensional Stiefel manifold of $p\times n$ matrix ${\bf U}_{1}$ with orthonormal columns
$\bf U_{1}^{T}{\bf U_{1}}={\bf I}$. Note that ${\bf V} \in {\mathcal O(p)}$, thus $\mbox{tr}{\mbSigma}^{-1}{\bf S}=
\mbox{tr}{\bf V}{\mbGamma}^{-1}{\bf V}^{\top}{\bf U}_{1}{\bf L}{\bf U}_{1}^{\top}=\mbox{tr}{\mbGamma}^{-1}{\bf H}_{1}
{\bf L}{\bf H}_{1}^{\top}$, where ${\bf H}_{1}={\bf V}^{\top}{\bf U}_{1} \in {\mathcal V_{n, p}}$. Note that ${\bf H}_{1}
\notin {\mathcal O(p)}$, the von Neumann Theorem then might possibily fail to be true in general.

\vspace{0.3cm}
\noindent {\it 6.2. When $p > n$, both $n$ and $p$ are large so that $c \in (1, \infty)$}
\vspace{0.2cm}

\indent When $p > n$, since the population covariance matrix $\mbSigma$ is assumed to be a positive definite symmetric
matrix, so its $p$ eigenvalues are all positive, however rank(${\bf S})=n$. Thus, it has $p-n$ sample eigenvalues being $0$
in probability. As such, it seems difficult to get all the consistent eigenvalue estimators of the population eigenvalues.
When the sample size $n$ and the dimension $p$ are all large, maybe we only need the $n$ largest eigenvalue estimators to be
consistent for the first $n$ largest population eigenvalues. If this is the case, then the method developed in this note is
still applicable.

\vspace{0.3cm}
\noindent {\it 6.3. Whither the high-dimensional low-saample size (HDLSS) categorical data models?}
\vspace{0.2cm}

\indent When $p >> n$, namely for the HDLSS categorical models, our method might still have of help in some situations.
In the context of HDLSS categorical models, abundant in genomics and bioinformatics, with relatively smaller sample size
$n$ but also often $p >> n$. Motivated by the 2002-03 severe acute respiratory syndrome coronavirus (SARSCoV) epidemic
model, a general model of comparing $G~(\geq 2)$ groups is considered by Sen, Tsai, and Jou (2007). Each sequence has $P$
positions, each one relating to a categorical response indexed as $1, \cdots, C$, and there are $n_g$ sequences in the
$g$th group, for $g = 1, \cdots, G$. For the $g$th group, $p$th position and $c$th category, let $n_{gpc}$ be the number
of sequences, and let $n_{gp}=\sum_{c=1}^{C} n_{gpc}$, for $p = 1, \cdots, P, g =1, \ldots, G$. Note that if there is no
missing value, each sequence, at each position, takes on one of the $C$ responses $1, \ldots, C$, so that $n_{gp} = n_g$,
for all $p =1,\cdots, P$. The combined group sample size is $n = \sum_{g=1}^G n_g$. For geographically separated sequences,
the assumption of independence of the $G$ groups could be reasonable, but the sequences within a group may not be
independent due to their shared ancestry. For SARSCoV or HIV genome sequences, because of rapid evolution of the virus,
the independence assumption may not be very stringent. Further, for each sequence, the responses at the $P$ positions are
generally not independent nor necessarily identically distributed.

\indent For SARSCoV genome sequences, the scientific focus is the statistical comparison of different strata to coordinate
plausible differences to pertinent environmental factors. In many fields of applications, particularly, in
genomic studies, not only do we have $p > > n$, but also often $n$ small, leading to a curse of dimensional problems. One
encounters conceptual and operational roadblocks due to too many unknown parameters. For such genomic sequences, any single
position (gene) yields very little statistical information. Hence, a composite measure of qualitative variation over the
entire sequence is sought to be a better way of gauzing statistical group discrimination. In the specific context, some
molecular epidemiologic studies have advocated suitable external sequence analysis like multivariate analysis of variance
(MANOVA), although there are impasses of various types. Genomic research is a prime illustration for motivating
appropriate statistical methodology for comprehending the genomic variation in such high dimensional categorical data models.
Variation (diversity) in such large $P$ small $n$ models can not be properly statistically studied by standard discrete
multivariate analysis tools, using the full likelihood approach. For qualitative data models, the Gini-Simpson (GS)
index (Gini, 1912; Simpson, 1949) and Shannon entropy (Shannon, 1948) are commonly used for statistical analysis in many
other fields, including genetic variation studies (Chakraborty and Rao, 1991). The Hamming distance provides an
average measure that does not ignore dependence or possible heterogeneity. The U-statistics methodology (Hoeffding, 1948)
is incorporated to obtain optimal nonparametric estimators and their jackknife variance estimators. The distribution theory
would have followed the general results of Tsai and Sen (2005). Because of underlying restraints, a pseudo-marginal
diversity function approach based on Hamming distance is considered by Sen, Tsai, and Jou (2007) in a statistical inference
setup. However, in the present context, $P$ is very large while the $n_g$ are all small. This requires the exploration of
asymptotics for the $P > > n$ environment that is considered in the following. In the genome sequence context, we are
confronted with the $P > > n$ environment. Within this framework, we encounter two scenarios: (i) $P > > n$; $n$ at least
moderately large, and (ii) $P > > n$ with $n$ small. In (i), the sample estimates of the Hamming distances are all
U-statistics, standard asymptotics (Sen 1977, Tsai and Sen 2005) hold: The estimators are asymptotically (as $n \to \infty$)
normal and their jackknifed variance estimators are consistent. Hence, we shall not enter into a detailed discussion of (i).
Case (ii), more commonly encountered in genomic studies, has different perspectives. We need to use the appropriate central
limit theory (CLT) for dependent sequences of bounded random variables (for the details see Sen, Tsai, and Jou, 2007).

\indent Tsai and Sen (2010) showed that the Shannon entropy is more informative than the GS index in the sense of the
Lorenz ordering makes it more appealing to consider the Shannon entropy. For HDLSS genomic models, they suspected that
the information might not be fully captured in a pseudo-marginal setup. The Hamming-Shannon pooled measures are more
informative than the pseudo-marginal Hamming-Shannon measures. To capture greater information, some new genuine
multivariate analogs of Shannon entropy are proposed.

\indent For HDLSS categorical models, Tsai and Sen (2010) showed that the Hamming-Shannon measures have the properties of
nested subset monotonicity and subgroup decomposability. Let $\pi_{gpc}$ denote the $c$th cell probability for the $p$th
marginal law ${\mbpi_{gp}}=(\pi_{gp1,}, \cdots, \pi_{gpC})^{'}$ of group $g ~(1 \leq c \leq C, 1 \leq g \leq G, 1 \leq p
\leq P)$, and let $n_{gpc}$ be the cell frequencies for the $p$th marginal table corresponding to the $g$th group, so that
the MLE of $\pi_{gpc}$ is $\widehat \pi_{gpc}=n_{gpc}/n_{g}, 1 \leq c \leq C$, where $n_g=\sum_{c=1}^{C}n_{gpc}$, the same
for every $p~(=1, \cdots, P)$. We incorporate the jackknife methodology to obtain the nonparametric estimators.
The jackknife estimator, the plug-in estimator based on the MLE of ${\mbpi_{gp}}$, of the Hamming-Shannon measure is
considered. The difficulties of the HDLSS asymptotics in the HDLSS genomic context are assessed and suitable permutation
procedures are appraised along with. Under the null hypothesis, the homogeneity of the $G$ groups, the advantage of the
resulting permutation invariance structure is taken. Therefore, we proceed with this extended permutation-jackknife
methodology.

\indent Consider all possible equally likely permutations of the observations for each $p$, each having the same conditional
probability $\frac{1}{N}$, where $N=n!/\prod_{g=1}^{G} n_{g}!, g=1,\cdots, G$. Let ${\bf Y}_{1}=(T_{2:1}, \cdots,
T_{G: G-1})^{t}$ with $T_{\{i: i-1\}}$ being defined in expression (63) (Tsai and Sen, 2010). In practice, to overcome the
difficulty that $N$ is too large we may choose $N_1$, which is sufficiently large but $N_1 < < N$, instead. Next, generate
a set of $(N_{1}-1)$ permutations. For this construction, we use the permutation distribution generated by the set of all
possible permutations among themselves. Consider ${\bf Y}_{i}$ be the $(i-1)$th corresponding permutation of ${\bf Y}_{1}$,
$i=2, \cdots, N_1$, and the corresponding covariance matrix ${\bf S}_{N_1}=({N_1}-1)^{-1}\sum_{i=1}^{N_1}({\bf Y}_{i}-
\bar {\bf Y})({\bf Y}_{i}-\bar {\bf Y})^{t}$, where $\bar {\bf Y}={N_1}^{-1}\sum_{i=1}^{N_1}{\bf Y}_{i}$. In practice,
$N_1$ is taken larger enough so that $lim_{N_1 \to \infty} P/N_{1}=c \in (0, 1)$. Do the spectral decomposition for the
matrix ${\bf S}_{N_1}$, then with the eigenvalues of ${\bf S}_{N_1}$ being replaced by the new corresponding eigenvalues
obtained based on equation (4.11), just the same as that of the sample covariance matrix ${\bf S}$ being replaced
by $\Hat{\mbSigma}_{T}$, to get the new improved Jackknife covariance matrix to replace ${\bf S}_{N_1}$ for the statistical
inference. Thus, under these circumstances,  our procedure proposed in Section 4 works well for the HDLSS categorical
data models.

\indent  In Section 4, when the sample size $n$ is larger than the dimension $p$ so that $\lim_{n \to \infty}p/n=c \in
(0, 1)$, we may note that $n\Hat{\mbSigma}_{T}$ is Wishart distributed with mean matrix ${\mbSigma}$. It is demonstrated
that $\Hat{\mbSigma}_{T}$ is the MLE of the population covariance matrix ${\mbSigma}$, hence, it is consistent. Moreover,
it is easy to see that $\Hat{\mbSigma}_{T}$ reduces to the sample covariance matrix ${\bf S}$ when the dimension $p$ is
fixed. Hence, when $n > P$ with the proposed novel estimator $\Hat{\mbSigma}_{T}$ replacing the sample covariance matrix
${\bf S}$, the traditional case of fixed dimension $p$ and the modern case of high-dimensional setup can then be integrated
into a unified theory. Thus, when $n > P$, the proposed novel estimator $\Hat{\mbSigma}_{T}$ of the population covariance
matrix ${\mbSigma}$ plays the fundamental role for further theoretical development of statistical inference. Practically,
it does have of help for some HDLSS categorical data models. When the sample size $n$ is moderate, Sen, Tsai and Jou (2007)
proposed the optimal nonparametric methods for the genomic data. When $P >> n$ and the sample $n$ is small, Tsai, and Sen
(2010) incorporated the permutation and Jackknife methodology to make statistical inferences for the genomic data.
Hopefully, the optimal statistical methods can be of help for scientific breakthroughs and also for real-world applications
of gene science.

\vspace{0.3cm}
\noindent {\bf Acknowledgements.}
\vspace{0.3cm}

\vspace{0.3cm}
\noindent {\bf References}
\vspace{0.3cm}

\begin{enumerate}
\item Anderson, T. W. (2003). {\it An Introduction to Multivariate Statistical Analysis,} 3rd edition,
                Wiely, New York.

%\item Bai, Z. D. and Silverstein, J. (2010). {\it Special Analysis of Large Dimensional Random Matrices},
%Springer, New York.

\item Bai, Z. D., Miao, B. Q., and Pan, G. H. (2007). On asymptotics of eigenvectors of large sample covariance matrix.
     Ann. Probab. 35, 1532-1572.

\item Chakraborty, R., and Rao, C. R. (1991). Measurement of genetic variation for evolutionary studies. In {\it Handbook
   of Statistics Vol. 8: Statistical Methods in Biological and Medical Sciences} (Eds., C. R. Rao and R. Chakraborty),
   pp. 271-316. North-Holland.

\item Choi, S. I. and Silverstein, J. W. (1995). Analysis of the limiting spectral distribution of large
               dimensional random matrices. J. Multivariate Anal. 54, 295-309.

\item Gini, C. W. (1912). Variabilita e mutabilita. Studi Economico-Giuridici della R. Universita de Cagliary 2, 3-159.

\item Hoeffding, W. (1948). A class of statistics with asymptotically normal distribution. Ann. Math. Statist. 19, 293-325.

\item James, W. and Stein, C. (1961). Estimation with quadratic loss. Proc. Fourth Berkeley Symp. Math.
               Statist. Probab. 1, 361-379. California Press, Berkeley, CA.

\item Johnstone, I. M. and Lu, A.Y. (2009). On consistency and sparsity for principal components analysis in
                 high dimensions. J. Amer. Statist. Assoc. 104, 682-693.

\item Johnstone, I. M. and Paul, D. (2018). PCA in high dimensions: An orientation. Proceedings of the IEEE.
                 106, 1277-1292.

\item Ledoit, O. and P$\acute{e}$ch$\acute{e}$, S. (2011). Eigenvectors of some large sample covariance matrix
                   ensembles. Probab. Theory Relat. Fields. 151, 233-264.

\item Ledoit, O. and Wolf, M. (2012). Nonlinear shrinkage estimation of  large-dimensional covariance matrices.
              Ann. Statist. 40, 1024-1060.

\item Ledoit, O. and Wolf, M. (2018). Optimal estimation of a large-dimensional covariance matrix under Stein's loss.
             Bernoulli 24, 3791-3832.

\item Mar$\check{c}$enko, V. A. and Pastur, L. A. (1967). Distribution of eigenvalues for some sets of random matrices.
                Sb. Math. 1, 457-483.

%\item Perlman, M. D. (2007). {\it Multivariate Statistical Analysis,} Univ. Washington, Seattle, Washington.

%\item Muirhead, R. J. (1982). {\it Aspects of Multivariate Statistical Theory}, Wiley, New York.

\item Rajaratnam, B. and Vincenzi, D. (2016). A theoretical study of Stein's covariance estimator. Biometrika 103, 653-666.

\item Pourahmadi, Mohsen. (2013). {\it High-Dimensional Covariance Estimation}, Wiley, New York.

\item Sen, P. K. (1977). Some invariance principles relating to Jackknifing and their role in sequential analysis.
             Ann. Statist. 5, 315-329.

\item Sen, P. K., Tsai, M.-T., and Jou., Y. S. (2007). High-dimension, Low-sample size perspectives in constrained
                statistical inference: The SARSCoV RNA genome in illustration. J. Amer. Statist. Assoc. 102, 686-694.

\item Shannon, C. E. (1948), A mathematical theory of communication, Bell System Techni. J. 27, 379-423, 623-656.

\item Simpson, E. H. (1949), The measurement of diversity. Nature, 163, 688.

\item Silverstein, J. W. (1995). Strong convergence of the empirical distribution of eigenvalues of large dimensional
               random matrices. J. Multivariate Anal. 55, 331-339.

\item Stein, C. (1956). Inadmissibility of the usual estimator of the mean of a multivariate normal distribution.
               Proc. Third Berkeley Symp. Math. Statist. Probab. 1, 197-206. California Press, Berkeley, CA.

\item Stein, C. (1975). Estimation of a covariance matrix. Rietz lecture, 39th Annual Meeting IMS, Atalanta, Georgia.

\item Stein, C. (1986). Lectures on the theory of estimation of many parameters. J. Math. Sci. 43, 1373-1403.

\item Tsai, M.-T. (2018). On the maximum likelihood estimator of a covariance matrix. Math. Method. Statist. 27, 71-82.

\item Tsai, M.-T and Sen, P. K. (2005). Asymptotically optimal tests for parametric functions against ordered functional
            alternatives. J. Multivariate Anal. 95, 37-49.

\item Tsai, M.-T. and Sen, P. K. (2010). Entropy based constrained inference for some HDLSS genomic models: UI tests in
            Chen-Stein perspective. J. Multivariate Anal. 101, 1559-1573.

\item Tsai, C.-H. and Tsai, M.-T. (2024a). On the decomposite $T^{2}$-test when the dimension is large. arXiv.2403.01516.

\item Tsai, M.-T. and Tsai, C.-H. (2024b). On the orthogonally equivariant estimators of a covariance matrix.
                arXiv.2405.06877.

\item Zagidullina, A. (2021). {\it High-Dimensional Covariance Matrix Estimation: An Introduction to Random Matrix Theory},
                   SpringerBriefs in Applied Statistics and Econometrics. Switzerland.

\item Uhlig, H. (1994). On singular Wishart and singular multivariate Beta distributions. Ann. Statist. 22, 395-405.

\end{enumerate}

\end{document}